\documentclass{article}

\usepackage{amsmath,amsthm,amssymb,amsfonts}
\usepackage{latexsym}
\def\pmb#1{\setbox0=\hbox{$#1$}%
\kern-.025em\copy0\kern-\wd0
\kern.05em\copy0\kern-\wd0
\kern-.025em\raise.0433em\box0}

\newfont{\field}{msbm10 at 11pt}
\newcommand{\R}{\mbox{\field \symbol{82}}}          

\newtheorem{theorem}{Theorem}[section]
\newtheorem{lemma}{Lemma}[section]

\newtheorem{remark}{Remark}[section]

\begin{document}
\title{\textit{Power Mean Curvature Flow in Lorentzian
    Manifolds}}
\author{Guanghan Li$^{1}$ and Isabel M.C.\ Salavessa$^{2}$}
\date{}
\protect\footnotetext{\!\!\!\!\!\!\!\!\!\!\!\!
$1$ Partially supported by CNSF Grants no. 10501011 and by
FCT through
the Plurianual of Project no. 1143 of CFIF.\\
$2$ Partially supported by FCT through
POCI/MAT/60671/2004 and the Plurianual of CFIF.\\[1mm]
{\footnotesize Guanghan Li\\
Centro de Fisica das Interac\c{c}\~{o}es Fundamentais,
Instituto Superior T\'{e}cnico, \\
Edif\'{\i}cio Ci\^{e}ncia, Piso 3,
Av.\ Rovisco Pais, 
1049-001 Lisboa, Portugal,~~~~~and\\[1mm]
School of Mathematics and Computer Science,
Hubei University,\\
Wuhan, 430062, P. R. China;~~
e-mail: liguanghan@163.com}\\[2mm]
{\footnotesize Isabel M.C. Salavessa\\
 Centro de F\'{\i}sica das Interac\c{c}\~{o}es Fundamentais,
Instituto Superior
T\'{e}cnico, \\
Edif\'{\i}cio Ci\^{e}ncia,
Piso 3, Av.\ Rovisco Pais, 1049-001 Lisboa, Portugal;\\
e-mail: isabel@cartan.ist.utl.pt}
}
\maketitle 
\noindent
{\small {\bf Abstract:}
 We study the motion of an $n$-dimensional
closed spacelike hypersurface in a Lorentzian manifold in the
direction of its past directed normal vector, where the speed
equals a positive power $p$ of the mean curvature. We prove that
for any $p\in (0, 1]$, the flow exists for all time when the Ricci
tensor of the ambient space is bounded from below on the set of
timelike unit vectors. Moreover, if we assume that all evolving
hypersurfaces stay in a precompact region, then the flow converges
to a stationary maximum spacelike hypersurface.}\\[2mm]
\markright{\sl\hfill Li -- Salavessa \hfill}

\section{Introduction}
\renewcommand{\thesection}{\arabic{section}}
\renewcommand{\theequation}{\thesection.\arabic{equation}}
\setcounter{equation}{0}

Curvature flow has been studied extensively for more than twenty
years. Mean curvature flow was first studied in a well-known paper
of Huisken \cite {h}. In \cite {a1}, Andrews studied the motion of
a closed convex hypersurface in $\R^{n+1}$ with the normal
velocity proportional to a positive power of Gaussian curvature.
While in \cite {s1}, Schulze considered the power mean curvature
flow of convex hypersurfaces also in the Euclidean space. In the
case of Lorentzian geometry, Ecker, Huisken \cite {eh}, and
Gerhardt \cite {g2, g3, g4, g5} studied several kinds of curvature
flow of spacelike hypersurfaces in order to find hypersurfaces
with prescribed curvatures.

In this paper, we shall consider the motion of a closed spacelike
hypersurface in a Lorentzian manifold, in the direction of its
past directed normal vector field, such that the evolution
velocity is a positive power of the mean curvature. Precisely, let
$N$ be a globally hyperbolic $(n+1)$-dimensional Lorentzian
manifold with a compact Cauchy hypersurface $S_0$, and $M$ an
$n$-dimensional smooth manifold. We shall investigate the
following problem. Let $M_0\subset N$ be a spacelike hypersurface
the mean curvature of which is strictly positive, we then consider
the power mean curvature flow
\begin{equation}
\dot x= sign(p)\, H^p\, \nu,
\end{equation}
with initial data $x(0)=x_0$, where $sign(p)$ is the sign function
of $p$, and $x_0$ an embedding of initial hypersurface
$M_0=x(0)(M)$.  Here $H$ is the mean curvature of the flow
hypersurfaces $M(t)=x_t(M)$, with respect to the past directed
normal $\nu$ and $p$ a non-zero real number. The evolution
equation (1.1) is a parabolic problem, the solution exists for a
short time. If we assume that $M_0$ is a graph over $S_0$, then we
can also write $M(t)$ as a graph over $S_0$, i.e.
$$M(t)=\mbox {graph}\, u(t,x),$$
where $x\in S_0$ is an abbreviation of the space-like components.

For $p=1$, Ecker, Huisken \cite {eh} and Gerhardt \cite {g2}
studied the evolution problem $\dot x= (H-f)\nu$, in order to find
a spacelike hypersurface with a given mean curvature function $f$
provided so-called upper and lower barriers for $(H, f)$ are
satisfied. If $N$ satisfies a mean curvature barrier condition and
strong volume decay condition, Gerhardt \cite {g5} proved that the
inverse mean curvature flow (the case $p=-1$ in (1.1)) exists for
all time and provides a foliation of the future of the initial
hypersurface $M_0$. Furthermore, if $N$ is a asymptotically
Robertson-Walker space, Gerhardt \cite {g6} showed that the leaves
of an inverse mean curvature flow converge to an umbilical
hypersurface and provide transition from big Crunch to big Bang.

If the Ricci tensor of $N$ is bounded from below on the set of
timelike unit vector, we shall prove that the power mean curvature
flow (1.1), with $p\in (0, 1]$ exists for all time, and converges
to a stationary maximum spacelike hypersurface when all evolving
hypersurfaces stay in a precompact region, i. e. we have the
following theorem.
\begin{theorem}
Let $N$ be a cosmological spacetime with compact Cauchy
hypersurface $S_0$ and $M_0$ a closed spacelike hypersurface with
positive mean curvature in $N$. If the Ricci tensor of $N$ is
bounded from below on the set of timelike unit vectors
 \begin{equation}
\overline R_{\alpha\beta}\nu^{\alpha}\nu^{\beta}\geq -\Lambda,
\quad ~~~\forall <\nu, \nu>=-1,
\end{equation}
for a positive constant $\Lambda$, then the power mean curvature
flow equation (1.1) with initial hypersurface $M_0$ exists for all
time, when $p\in(0,1]$. Moreover if we further assume that the
limit $\lim_{t\rightarrow \infty}\inf_{x\in S_0}u(t,x)$ is bounded
from below,  then the solutions to (1.1) converge to a stationary
maximum spacelike hypersurface as $t\rightarrow \infty$.
\end{theorem}
\begin{remark} $(i)$ That the
limit $\lim_{t\rightarrow \infty}\inf_{x\in S_0}u(t,x)$ is bounded
from below is equivalent to say that all evolving hypersurfaces
stay in a precompact region of $N$. Moreover in this case,  the
flow converges to stationary maximum spacelike hypersurface
without the assumption on the Ricci tensor bound (1.2) (see Lemma
4.1).
$(ii)$ If we suppose that the upper and lower barriers for the mean
curvature hold like that in \cite {g2}, the results in Theorem 1.1
still hold.
\end{remark}
For the notations and definitions we refer to section 2 for more
detailed description. In section 3 we look at the curvature flow
associated with our problem, and corresponding evolution equations
for the basic geometric quantities of the flow hypersurfaces. In
section 4 lower bound estimates for the evolution problem are
proved, while a priori estimate for the norm of the second
fundamental form is derived in section 5. Finally in section 6, we
demonstrate that the maximum existence time is infinity and the
evolution equation converges to stationary solution, which
completes the proof of Theorem 1.1.

\section{Preliminaries}
\renewcommand{\thesection}{\arabic{section}}
\renewcommand{\theequation}{\thesection.\arabic{equation}}
\setcounter{equation}{0}

In this section we state some basic concepts and formulas which
can be found in \cite {g2} or \cite {g3}. Let $N$ be a globally
hyperbolic Lorentzian manifold with a compact Cauchy surface
$S_0$. $N$ is a topological product $\mathbb{R}\times S_0$, where
$S_0$ is a compact Riemannian manifold, and there exists a
Gaussian coordinate system $\{x^{\alpha }\}$, and such that the
metric in $N$ has the form
\begin{equation}
d\overline s_N^2=e^{2\psi}\left({-dx^0}^2+\sigma _{ij}(x^0,
x)dx^idx^j\right),
\end{equation}
where $\sigma _{ij}$ is a Riemannian metric, $\psi$ a function on
$N$, and $x$ an abbreviation for the space-like components
$(x^i)$. We also assume that the coordinate system is future
oriented, i. e. the time coordinate $x^0$ increases on future
directed curves.

Let $M$ be a space-like hypersurface, i. e. the induced metric is
Riemannian, with a differentiable normal $\nu $ that is timelike.
Let $(\xi ^i)$ be the local coordinates on $M$. Geometric
quantities on $N$ will be denoted with a bar, for
example, by $(\overline g_{\alpha
\beta})$, $(\overline R_{\alpha \beta \gamma \delta})$, etc., and
those in $M$ without the bar, as
 $(g_{ij})$, $(R_{ijkl})$, etc. Greek indices range
from $0$ to $n$ and Latin from $1$ to $n$. For a function $u$ in
$N$, $(u_{\alpha })$  and $(u_{\alpha \beta
})$ denotes respectively the gradient and the Hessian of $u$.
We use the summation index convention on repeated indices.
We assume the above time-like normal $\nu$ is past directed. Then
the Gaussian formula of $M$ is given by
\begin{equation}
x_{ij}^{\alpha}=h_{ij}\nu ^{\alpha}.
\end{equation}
We represent covariant derivative of a tensor as
a full  tensor. For instance
\begin{equation}
x_{ij}^{\alpha}=x_{,ij}^{\alpha}-\Gamma
_{ij}^kx_k^{\alpha}+\overline {\Gamma}_{\beta
\gamma}^{\alpha}x_i^{\beta}x_j^{\gamma},
\end{equation}
where the comma indicates ordinary partial derivatives, and
$\{\Gamma _{ij}^k\}$, $\{\overline \Gamma _{\beta
\gamma}^{\alpha}\}$ the Christoffel symbols of $M$ and $N$,
respectively. The quantities
$\overline R_{\alpha \beta \gamma
\delta;\varepsilon}$ denotes the covariant derivative of the curvature
tensor, where $i=\frac{\partial}{\partial \xi^i}=x^{\varepsilon}_i
\varepsilon=
x^{\varepsilon}_i\frac{\partial}{\partial x^{\varepsilon}}$.
The  definition of the second fundamental form $(h_{ij})$
is taken with respect to $\nu$. The Weingarten equation is given
by
\begin{equation}
\nu _i^{\alpha}=h_i^kx_k^{\alpha},
\end{equation}
where
$\nu_i^{\alpha}=\nu _{,i}^{\alpha}+\nu ^{\gamma}\overline
{\Gamma}_{i \gamma}^{\alpha}$. Recall the Codazzi and the
Gaussian equations,
\begin{eqnarray}
&&h_{ij,k}-h_{ik,j} =\overline R_{\alpha \beta \gamma \delta} \nu
^{\alpha}x_i^{\beta}x_j^{\gamma}x_k^{\delta}\\
&&R_{ijkl}=-\left(h_{ik}h_{jl}-h_{il}h_{jk}\right)
+\overline R_{\alpha \beta
\gamma \delta}x_i^{\alpha}x_j^{\beta}x_k^{\gamma}x_l^{\delta}.
\end{eqnarray}
Let $M=\mbox {graph}u|_{S_0}$ be a space-like hypersurface
\begin{equation*}
M=\{(x^0, x)|x^0=u(x), x\in S_0\},
\end{equation*}
The induced metric on $M$ has the form
\begin{equation*}
g_{ij}=e^{2\psi}\left(-u_iu_j+\sigma _{ij}\right),
\end{equation*}
where $\sigma _{ij}$ are fuctions of  $(u, x)$, and $u_i=\frac
{\partial u}{\partial x^i}$. The inverse  matrices
$(g^{ij})=(g_{ij})^{-1}$ and $(\sigma ^{ij})=(\sigma _{ij})^{-1}$
are related by
\begin{equation*}
g^{ij}=e^{-2\psi}\left(\frac {u^i}{v}\frac {u^j}{v}+\sigma ^{ij}\right),
\end{equation*}
where
\begin{equation*}
u^i=\sigma ^{ij}u_j, \qquad v^2=1-\sigma
^{ij}u_iu_j\equiv1-|Du|^2.
\end{equation*}
Hence, graph$u$ is space-like if and only if $|Du|<1$.
The contravariant past directed normal vector and its
covariant form are respectively given by
\begin{eqnarray}
(\nu ^{\alpha})&=&-v^{-1}e^{-\psi}(1, u^i)=-v^{-1}e^{-\psi}(1, u^1,
\cdots, u^n)\\
(\nu _{\alpha})&=&v^{-1}e^{\psi}(1, -u_i)=v^{-1}e^{\psi}(1, -u_1,
\cdots, -u_n).\nonumber
\end{eqnarray}
From (2.3) and (2.7) one obtains
\begin{equation}
e^{-\psi}v^{-1}h_{ij}=-u_{ij}-\overline {\Gamma
}_{00}^0u_iu_j-\overline {\Gamma }_{0j}^0u_i-\overline {\Gamma
}_{0i}^0u_j-\overline {\Gamma }_{ij}^0.
\end{equation}
Set $(\overline h_{ij})$ be the second fundamental form of the
hypersurface $\{x^0=const.\}$. Then
\begin{equation}
-\overline {\Gamma}_{ij}^0=e^{-\psi}\overline h_{ij},
\end{equation}
and by computing the Christoffel symbol
of $N$ one derives
\begin{equation}
\overline h_{ij}e^{-\psi}=-\frac 12 \dot{\sigma} _{ij}-\dot
{\psi}\sigma _{ij},
\end{equation}
where the dot indicates differentiation with respect to $x^0$.

In \cite{g2,g3,g4} it is defined a Riemannian metric $(\widetilde{g}_{\alpha
\beta })$ by
\begin{equation*}
\widetilde{g}_{\alpha \beta
}dx^{\alpha}dx^{\beta}=e^{2\psi}\left({dx^0}^2+\sigma _{ij}dx^idx^j\right)
\end{equation*}
and corresponding norm of a vector field $\eta$ by
\begin{equation*}
|||\eta |||=(\widetilde{g}_{\alpha \beta }\eta ^{\alpha}\eta
^{\beta})^{\frac 12}
\end{equation*}
with similar notations for higher tensors.
\section{Evolution Equations}
\renewcommand{\thesection}{\arabic{section}}
\renewcommand{\theequation}{\thesection.\arabic{equation}}
\setcounter{equation}{0}

From now on, we always assume $p>0$, and so $sign(p)=1$. In order
to study the evolution problem (1.1), it is convenient to consider
the evolution equation
\begin{equation}
\dot x= (H^p-\tau )\nu,
\end{equation}
where $\tau $ is a small positive number. The evolution problem
(3.1) is a parabolic problem, hence a solution exists on a maximum
time interval $[0, T^*)$, $0<T^*\leq \infty$. In the following we
show how the metric, the second fundamental form and the normal
vector of the hypersurfaces $M(t)$ evolve. All time derivatives
are total derivatives. Here is just the special case, we refer to
\cite {g3} for more general results, so we omit the proofs of the
following lemmas.
\begin{lemma}
The metric, volume element, the normal vector, and the second
fundamental form of $M(t)$ satisfy the evolution equations
\begin{eqnarray*}
\dot{g}_{ij} & = & 2(H^p-\tau )h_{ij},\\
\dot{\nu} & = & \nabla_M(H^p)=g^{ij}(H^p)_ix_j\\
\dot{h}_i^j & = & (H^p)_i^j-(H^p-\tau) h_i^kh_k^j-(H^p-\tau)
\overline
          R_{\alpha \beta \gamma \delta}\nu ^{\alpha}x_i^{\beta}\nu
              ^{\gamma}x_k^{\delta}g^{kj},\\
\dot{h}_{ij} & = &
          (H^p)_{ij}+(H^p-\tau )h_i^kh_{kj}-(H^p-\tau)\overline
          R_{\alpha \beta \gamma \delta}\nu ^{\alpha}x_i^{\beta}\nu
              ^{\gamma}x_j^{\delta}.
\end{eqnarray*}
\end{lemma}
\noindent
Let $\|A\|^2=h_{ij}h^{ij}$.
\begin{lemma}
The mean curvature $H$ evolves according to the following
equations
\begin{eqnarray}
\frac d{dt}H &=& pH^{p-1}\Delta
H+p(p-1)H^{p-2}||\nabla_MH||^2\nonumber \\[-1mm]
&&-(H^p-\tau )(||A||^2+\overline
R_{\alpha \beta }\nu ^{\alpha}\nu^{\beta})\\
\frac d{dt}H^p &=& pH^{p-1}\Delta
H^p-pH^{p-1}(||A||^2+\overline
R_{\alpha \beta }\nu ^{\alpha}\nu
^{\beta})(H^p-\tau ).~~~
\end{eqnarray}
\end{lemma}
\begin{lemma}
The mixed tensor $h_i^j$ satisfies the parabolic equation
\begin{eqnarray}
\frac d{dt}h_i^j& = & pH^{p-1}\Delta
h_i^j-pH^{p-1}(||A||^2+\overline
R_{\alpha \beta }\nu ^{\alpha}\nu
^{\beta})h_i^j \nonumber\\
&&+(p-1)H^ph_i^kh_k^j
+\tau h_i^kh_k^j
+  p(p-1)H^{p-2}\nabla _iH\nabla ^jH\nonumber\\
&&+2pH^{p-1}\overline
R_{\alpha \beta \gamma
\delta}x_m^{\alpha}x_i^{\beta}x_k^{\gamma}x_r^{\delta}h^{km}g^{rj}
+\tau \overline R_{\alpha \beta \gamma
\delta}\nu ^{\alpha}x_i^{\beta}\nu ^{\gamma}x_m^{\delta}g^{mj}
\nonumber\\
&&- pH^{p-1}g^{kl}\overline
R_{\alpha \beta\gamma \delta}x_m^{\alpha}x_k^{\beta}x_l^{\delta}
(x_r^{\gamma}h_i^mg^{rj}+x_i^{\gamma}h^{mj})\nonumber\\
&&+(p-1)H^p\overline R_{\alpha \beta \gamma\delta}
\nu ^{\alpha}x_i^{\beta}\nu ^{\gamma}x_m^{\delta}g^{mj}\nonumber\\
&&+  pH^{p-1}g^{kl}\overline R_{\alpha \beta \gamma
\delta;\varepsilon}(\nu^{\alpha}x_k^{\beta}
x_l^{\gamma}x_i^{\delta}x_m^{\varepsilon}g^{mj}+
\nu^{\alpha}x_i^{\beta}
x_k^{\gamma}x_m^{\delta}x_l^{\varepsilon}g^{mj}).~~~~~~~~~
\end{eqnarray}
\end{lemma}
\noindent
We immediately deduce from (3.3)
\begin{lemma}
 If $H^p\geq \tau$ at $t=0$, then for any $t\in [0, T^*)$,
$H^p(t)\geq \tau$.
\end{lemma}

\section{Lower Bound Estimate}
\renewcommand{\thesection}{\arabic{section}}
\renewcommand{\theequation}{\thesection.\arabic{equation}}
\setcounter{equation}{0}

 The evolution
problem (3.1) exists on a maximum time interval $I=[0, T^*)$. We
will prove that $T^*=\infty$.
Because of the short time existence, and the initial hypersurface
$M_0$ is a graph over $S_0$, we can write
\begin{equation*}
M(t)=\mbox {graph}\, u(t)
=\{\left( u(t, x(t)),x(t)\right): ~x(t)\in S_0\},
\qquad \forall t\in I,
\end{equation*}
where $u$ is defined in the cylinder $Q_{T^*}=I\times S_0$.
From (3.1) and
using (2.7) for $\alpha =0$ one sees
that $u$ satisfies a parabolic equation of the form
\begin{equation}
\dot u = -e^{-\psi }v^{-1}(H^p-\tau),
\end{equation}
where $\dot u$ is a total derivative, i. e.
\begin{equation}
\dot u = \frac {\partial u}{\partial t}+u_i\dot x^i,
\end{equation}
and so
\begin{equation}
\frac {\partial u}{\partial t} = -e^{-\psi }v(H^p-\tau).
\end{equation}
Consequently, from Lemma 3.4,
 $\frac {\partial u}{\partial t}$ is non positive.
 Next we shall prove that the flow stays in a precompact
region in finite time.
\begin{lemma}
Let $N$ be a cosmological spacetime with a compact Cauchy
hypersurface and satisfy condition (1.2). Then, for any finite
$T$, $0<T\leq T^*$, the flow (3.1) stays in a precompact region
$\Omega_T$ for $0\leq t\leq T$.
\end{lemma}
\begin{proof}
First we claim that, we can choose a new time function
$\widetilde{x}^0$ such that the Lorentzian metric of $N$ also has
the form (2.1), and the conformal factor satisfies
$\widetilde{\psi}\geq 0$. Suppose that
$$\widetilde{x}^0=\zeta (x^0),\quad  \widetilde{x}^i=x^i, \quad i=1, 2,
\cdots, n$$ where $\zeta$ is a $C^1$ function with non-vanishing
derivative. Then the metric is given by in the new coordinates
$$d\overline s_N^2=e^{2\psi}(\dot{\zeta})^{-2}\left(
-(d\widetilde{x}^0)^2+
(\dot{\zeta})^{2}\sigma _{ij}(\zeta^{-1} (\widetilde{x}^0),
\widetilde{x})d\widetilde{x}^id\widetilde{x}^j\right), $$ where the dot
indicates the differentiation with respect to $\widetilde{x}^0$.
If we let
$$(\dot{\zeta})^{-2}=exp(-2\inf_{x\in S_0}\psi (x^0, x)),$$
and we then have $\widetilde{\psi} =\psi -\inf_{x\in S_0}\psi
(x^0, x))$, which is non-negative. This proves the claim, and
therefore without loss of generality, we may assume $\psi \geq 0$
in the metric (2.1).

From (4.3) and Lemma 3.4 $u$
is decreasing. Thus we only need to prove
that $u$ has a lower bound.
Set $\varphi
(t)=\inf _{S_0}u(x,t)$. We may assume there exists some $t_0\in (0,
T)$ such that $\varphi (t_0)<0$, otherwise $u$ were bounded, and
the lemma holds. The function $\varphi $ is Lipschitz continuous and
if $x_t$ is such that
the infimum is attained at $x_t$, then
$\frac {\partial \varphi (t)}{\partial t}=
\frac {\partial }{\partial t}u(t, x_t)$ holds for a.e. $t$,
(cf. \cite {g5}, Lemma 3.2).
 By (3.2) and (1.2), and using again Lemma 3.2 of \cite {g5}, we see that
$$\frac {d}{dt}\sup_{S_0}H\leq \Lambda (\sup _{S_0}H^p-\tau )
\leq \Lambda \sup _{S_0}H^p , $$
 which implies that
\begin{equation}
\sup
H^{1-p}(t)\leq \sup H^{1-p}(0)+(1-p)\Lambda t
\end{equation}
 for $p\in (0,1)$, and
 \begin{equation}\sup H(t)\leq
\sup H(0)e^{\Lambda t}
\end{equation} for $p=1$.
 We therefore have
$$\frac {\partial\varphi (t)}{\partial t}\geq -(\sup_{S_0}H^{1-p}(0)
+(1-p)\Lambda t)^{\frac p{1-p}}-\tau,
 \quad {\mbox for}\quad 0<p<1,$$
and
$$\frac {\partial \varphi (t)}{\partial t}\geq -C_1e^{\Lambda
t}\quad  {\mbox for}\quad  p=1.$$ From these inequalities we
immediately deduce for $0<p<1$
$$\varphi (t)\geq -\tau T-C_2
-C_1(\,\sup_{S_0} H^{1-p}(0)+(1-p)\Lambda T\,)^{\frac 1{1-p}},\quad
\forall 0\leq t\leq T,$$ and for $p=1$
$$\varphi (t)\geq \varphi (t_0)-C_1e^{\Lambda t},\quad
\forall 0\leq t\leq T,$$ proving the lemma.
\end{proof}

\section{$C^1$ Estimates}
\renewcommand{\thesection}{\arabic{section}}
\renewcommand{\theequation}{\thesection.\arabic{equation}}
\setcounter{equation}{0}

We consider a smooth solution of the evolution equation (3.1) in a
maximum time interval $[0, T^*)$. In order to prove that the
hypersurfaces remain uniformly space-like, we only have to prove that
the term
\begin{equation*}
\widetilde{v}=v^{-1}=\frac {1}{\sqrt{1-|Du|^2}}
\end{equation*}
is uniformly bounded in finite time, i.e. in $Q_T=[0,T]\times S_0$
for any $0<T<T^*$. We shall apply the maximum principle to the
evolution equation of the quantity $w=\widetilde{v}\varphi$, where
$\varphi $ is defined in (5.8) below. This was first used by
Gerhardt and the following proof is a slight modification in \cite
{g2}. Let $\eta$ be the covariant vector field $(\eta
_{\alpha})=e^{\psi}(-1, 0, \cdots, 0)$.
\begin{lemma}
The quantity $\widetilde{v}$ satisfies the the evolutive equation
\begin{eqnarray}
\lefteqn{\dot {\widetilde{v}}-pH^{p-1}\Delta \widetilde{v}=}\nonumber\\
 & = & -pH^{p-1}||A||^2\widetilde{v}-2pH^{p-1}h^{ij}
x_i^{\alpha}x_j^{\beta}\eta _{\alpha
           \beta}-pH^{p-1}g^{ij}\eta _{\alpha \beta
          \gamma}x_i^{\beta}x_j^{\gamma}\nu
          ^{\alpha}\nonumber\\
& & -pH^{p-1}\overline R_{\alpha \beta}\nu
          ^{\alpha}x_k^{\beta}\eta _{\gamma}x_l^{\gamma}g^{kl}
           -(p-1)H^p\eta _{\alpha\beta}\nu ^{\alpha }
\nu ^{\beta}-\tau \eta _{\alpha \beta}\nu ^{\alpha}\nu^{\beta}.~~~~~~
\end{eqnarray}
\end{lemma}
\begin{proof}
Covariant differentiation of
$\widetilde{v}=<\eta ,\nu >$ gives
\begin{eqnarray}
\widetilde{v}_i&=&\eta _{\alpha \beta}x_i^{\beta}\nu ^{\alpha}+\eta
_{\alpha}\nu _i^{\alpha},\nonumber\\
\widetilde{v}_{ij}&=&\eta _{\alpha \beta \gamma
}x_i^{\beta}x_j^{\gamma}\nu ^{\alpha}+\eta _{\alpha
\beta}x_{ij}^{\beta}\nu ^{\alpha}+\eta _{\alpha \beta}x_i^{\beta
}\nu _j^{\alpha}+\eta _{\alpha \beta}x_j^{\beta}\nu
_i^{\alpha}+\eta _{\alpha}\nu _{ij}^{\alpha},~~
\end{eqnarray}
and the time derivative of $\widetilde{v}$ is given by
\begin{eqnarray}
\dot {\widetilde{v}} & = & \eta _{\alpha \beta}\dot {
              x}^{\beta}\nu ^{\alpha}+\eta _{\alpha}\dot {\nu} ^{\alpha}=\eta
                _{\alpha \beta}\nu ^{\alpha }\nu ^{\beta}(H^p-\tau)+\eta _{\alpha
            }\nabla ^k(H^p)x_k^{\alpha} \nonumber\\
& = & \eta _{\alpha \beta }\nu ^{\alpha}\nu
             ^{\beta}(H^p-\tau)+pH^{p-1}\nabla ^kHx _k^{\alpha}\eta _{\alpha}.
\end{eqnarray}
Now by Weingarten formula (2.4) and using Codazzi equation (2.5)
we have
\begin{eqnarray}
\eta _{\alpha}\nu _{ij}^{\alpha}g^{ij} & = & \eta _{\alpha
             }g^{ij}g^{kl}(h_{il,j}x_k^{\alpha}+h_{il}x_{kj}^{\alpha})\nonumber\\
& = & \eta _{\alpha
             }g^{ij}g^{kl}\left(
x_k^{\alpha}[h_{ij,l}+\overline R_{\varepsilon \beta \gamma \delta }
             \nu ^{\varepsilon}x_i^{\beta}x_l^{\gamma}x_j^{\delta}]+h_{il}h_{kj}\nu
             ^{\alpha}\right) \nonumber\\
& = & \eta _{\alpha
             }\nabla ^kHx_k^{\alpha}+ \overline R_{\varepsilon \gamma }\nu
             ^{\varepsilon}\eta ^{\alpha}x_k^{\alpha}x_l^{\gamma}g^{kl}+\eta _{\alpha}\nu
             ^{\alpha}||A||^2.
\end{eqnarray}
Substituting (5.2), (5.3) and (5.4) into (5.1),  we get the lemma.
\end{proof}
\noindent
The following lemma is a result in \cite {g2}
\begin{lemma} Consider the flow in a precompact region $\Omega _T$.
 Then there is a constant
$c=c(\Omega_T)>0$ such that for any
positive function $0<\epsilon=\epsilon (x)$ on $S_0$ and any
hypersurface $M(t)$ of the flow we have
\begin{equation*}
|||\nu |||\leq c\widetilde{v},
\end{equation*}
\begin{equation*}
g^{ij}\leq c\widetilde{v}^2\sigma^{ij} \mbox{~~(as~ operators)},
\end{equation*}
\begin{equation*}
|h^{ij}\eta _{\alpha \beta}x_i^{\alpha}x_j^{\beta}|\leq \frac
{\epsilon} 2||A||^2\widetilde{v}+\frac {c}{2\epsilon}\widetilde{v}^3.\\[1mm]
\end{equation*}
\end{lemma}
\noindent Using the previous lemmas, and that
 $|H^p\eta _{\alpha \beta}\nu
^{\alpha}\nu^{\beta}|\leq cH^p\widetilde{v}^2$, $|||\nu |||\leq
c\widetilde{v} $, and the  use of the Young inequality and the
relation $H^2\leq n||A||^2$ leads to:
\begin{lemma}There is constant $c=c(\Omega_T)$ such that for any
positive function $0<\epsilon=\epsilon (x)$ on $S_0$, the term
$\widetilde{v}$ satisfies a parabolic inequality of the form
\begin{equation}
\dot {\widetilde{v}}-pH^{p-1}\Delta \widetilde{v} \leq
           -pH^{p-1}(1-\epsilon)||A||^2\widetilde{v}+pH^{p-1}\left(1+\frac
           {1}{\epsilon}\right)\widetilde{v}^3+c\widetilde{v}^2,
\end{equation}
\end{lemma}
\noindent
We stress that this constant $c$ depends on $\Omega_T$ and
 on all geometric quantities of the ambient space
restricted to $\Omega$ we have ben considering.
Now set $u^i=g^{ij}u_j$.
\begin{lemma}
Let $M(t)=\mbox {graph}\, u(t)$ be the flow hypersurfaces, then we
have
\begin{eqnarray}
\dot {u}-pH^{p-1}\Delta u &=& (p-1)H^pe^{-\psi}\widetilde{v}+\tau
e^{-\psi}\widetilde{v}\nonumber\\
&&+pH^{p-1}(-e^{-\psi}g^{ij}\overline
h_{ij}+\overline {\Gamma}_{00}^0||Du||^2+2\overline
{\Gamma}_{0i}^0u^i).~~~~~
\end{eqnarray}
\end{lemma}
\noindent
This can be easily derived using (2.8).
  The following is  a lemma in
\cite {g2}
\begin{lemma}
Let $M\subset \Omega$ be a graph over $S_0$, then
\begin{equation*}
|\widetilde{v}_iu^i|\leq c\widetilde{v}^3+||A||e^{\psi}||Du||^2.
\\[1mm]
\end{equation*}
\end{lemma}
\noindent
Finally we  obtain a uniform bound of $\widetilde{v}$,
the proof of which is an adaptation to our case of the proof
of Proposition 3.7 of \cite {g2}.
\begin{theorem}Let $\Omega \in N$ be precompact. Then as long as
the flow stays in $\Omega$, the term $\widetilde{v}$ remains
uniformly bounded
\begin{equation}
\widetilde{v}\leq c=c(\Omega, M_0).
\end{equation}
\end{theorem}
\begin{proof}
Let $\mu$, $\lambda$ be positive constants, where $\mu$ is
 small and $\lambda$ large, and take
\begin{equation}
\varphi=e^{\mu e^{\lambda u}}.
\end{equation}
We may assume that $u\geq 1$,
otherwise we replace in (5.8) $u$ by $u+c$ for some $c$ large
enough. We also assume that $\widetilde{v}\geq 1$. We will see that
choosing $\lambda $
and $\mu$ conveniently, then
 $w=\widetilde{v}\varphi $ will be  uniformly bounded.
From
$$\dot {\varphi}= \varphi \mu \lambda e^{\lambda u}\dot u,$$
$$\Delta \varphi = \varphi (\mu \lambda e^{\lambda
u})^2||Du||^2+\mu \lambda ^2\varphi e^{\lambda u}||Du||^2+\mu
\lambda \varphi e^{\lambda u}\Delta u,$$
and Lemma 5.2 and 5.4 we conclude that
$$
\dot {\varphi}-pH^{p-1}\Delta \varphi\leq c\mu \lambda
e^{\lambda u}(pH^{p-1})\widetilde{v}^2\varphi-pH^{p-1}\mu \lambda
^2e^{\lambda u}\left(1+\mu e^{\lambda u}\right)||Du||^2\varphi,
$$
where we use $0<p\leq 1$. This inequality, Lemma
5.3 and 5.5 imply
\begin{eqnarray*}
\lefteqn{\dot {w}-pH^{p-1}\Delta w  = (\dot
              {\widetilde{v}}-pH^{p-1}\Delta
              \widetilde{v})\varphi +\widetilde{v}(\dot {\varphi}-pH^{p-1}
\Delta \varphi
              )-2pH^{p-1}D\widetilde{v}D\varphi} \\
  & \leq & -pH^{p-1}(1-\epsilon)||A||^2\widetilde{v}\varphi
  +pH^{p-1}c\left(1+\frac 1{\epsilon}\right)
\widetilde{v}^3\varphi+pH^{p-1}c\mu \lambda e^{\lambda
        u}\widetilde{v}^3\varphi\\
& & -pH^{p-1}\mu \lambda ^2e^{\lambda u}
 \left(1+\mu e^{\lambda u}\right)
||Du||^2\widetilde{v}\varphi-2pH^{p-1}\mu \lambda e^{\lambda
u}\varphi \widetilde{v}_i u_j g^{ij}.
\end{eqnarray*}
Here and in the above computations, we have used that $H$ is
bounded in finite time by Lemma 3.4, (4.4) and (4.5).
Note that the last term is $\leq pH^{p-1}\mu \lambda e^{\lambda
u}\varphi \left(c\widetilde{v}^3+||A||e^{\varphi }||Du||^2\right)$. From
the above we lastly arrive at
\begin{eqnarray*}
\lefteqn{\dot {w}-pH^{p-1}\Delta w  \leq
-pH^{p-1}(1-\epsilon)||A||^2\widetilde{v}\varphi
+pH^{p-1}c\left(1+\frac 1{\epsilon}\right)\widetilde{v}^3\varphi}\\
& & +pH^{p-1}c\mu \lambda e^{\lambda
u}\widetilde{v}^3\varphi-pH^{p-1}\mu \lambda ^2e^{\lambda u}
 \{1+\mu e^{\lambda u}\}
||Du||^2\widetilde{v}\varphi \\
& &  +2pH^{p-1}\mu \lambda e^{\lambda
u}||A||e^{\psi}||Du||^2\varphi.
\end{eqnarray*}
We estimate the last term on the right hand side by
\begin{eqnarray*}
\lefteqn{2pH^{p-1}\mu \lambda e^{\lambda
u}||A||e^{\psi}||Du||^2\varphi\leq}\\
&\leq&
pH^{p-1}||A||^2\widetilde{v}\varphi(1-\epsilon)+\frac{pH^{p-1}}{1-\epsilon}\mu
^2\lambda ^2 e^{2\lambda
u}\widetilde{v}^{-1}e^{2\psi}||Du||^4\varphi,
\end{eqnarray*}
and since $e^{2\psi}||Du||^2\leq \widetilde{v}^2$, we conclude
\begin{eqnarray}
\dot {w}-pH^{p-1}\Delta w  &\leq &
              pH^{p-1}c\left(1+\frac
              1{\epsilon}\right)\widetilde{v}^3
\varphi+pH^{p-1}c\mu \lambda e^{\lambda
u}\widetilde{v}^3\varphi\nonumber\\
&   &+ \left(\frac 1{1-\epsilon}-1\right)pH^{p-1}\mu ^2\lambda
^2e^{2\lambda
u}\widetilde{v}||Du||^2\varphi \nonumber\\
&&-pH^{p-1}\mu \lambda ^2e^{\lambda
u}||Du||^2\widetilde{v}\varphi.~~~~~~~~~~~
\end{eqnarray}
By setting $\epsilon =e^{-\lambda u}$, from (5.9) we derive
\begin{eqnarray}
\dot {w}-pH^{p-1}\Delta w &\leq&
              pH^{p-1}c\mu \lambda e^{\lambda
u}\widetilde{v}^3\varphi+pH^{p-1}ce^{\lambda u}\widetilde{v}^3\varphi
\nonumber\\
& & +pH^{p-1}\left(\frac
    {\mu}{1-\epsilon}-1\right)\mu\lambda ^2e^{\lambda
    u}|Du||^2\widetilde{v}\varphi.
\end{eqnarray}
Choose $\mu =\frac 12$ and $\lambda _0$ so large that
$
\frac
    {\mu}{1-\epsilon}\leq \frac 34, \quad \forall \lambda \geq
    \lambda _0$.
Since $||Du||^2=(\widetilde{v}^2-1)e^{-2\psi}\geq
\widetilde{v}^2e^{-2\psi}\geq c\widetilde{v}^2$
 we conclude that
the last term on the right hand side of (5.10) is
less than
\begin{equation*}
-\frac 18\lambda ^2e^{\lambda u}||Du||^2\widetilde{v}\varphi
pH^{p-1}\leq -c\lambda ^2e^{\lambda u}\widetilde{v}^3\varphi
pH^{p-1}.
\end{equation*}
Now (5.10) reduces to the form
\begin{equation*}
\dot {w}-pH^{p-1}\Delta w \leq
              cpH^{p-1}(1+\lambda -\lambda ^2)e^{\lambda
              u}\widetilde{v}^3\varphi.
\end{equation*}
Choosing $\lambda $ large enough such that $\lambda >\lambda
_0$, and applying the parabolic maximum principle  to $w$ gives
$$w\leq c=c(|w(0)|_{S_0}, \lambda _0, \Omega).$$
This completes the proof of Theorem 5.1.
\end{proof}

\section{Higher Order Estimates and Proof of Theorem 1.1}
\renewcommand{\thesection}{\arabic{section}}
\renewcommand{\theequation}{\thesection.\arabic{equation}}
\setcounter{equation}{0}

In this section, we  prove that as long as the flow stays
in a precompact set $\Omega \subset N$, the norm of the second
fundamental form of the evolutive hypersurface is a prior bounded
by a constant depending only on $\Omega $ and the initial
hypersurface $M_0$.
\begin{lemma} Let $\Omega \subset N$ be a precompact region and
assume that the flow (3.1) stays in $\Omega $ during the
evolution. Then the principal curvatures of the evolution
hypersurfaces $M(t)$ are uniformly bounded.
\end{lemma}
\begin{proof} We follow close \cite {g2} to prove that
$\varphi =\sup \{h_{ij}\eta ^i\eta ^j:||\eta ||=1\}$
is uniformly bounded.
The principal curvatures are bounded from above, for $H$ is positive
(see Lemma 3.4).
Let $0<T<T^*$, and $x_0=x_0(t_0)$, with $0<t_0\leq T$, be a point
in $M(t_0)$ such that
\begin{equation*}
\sup_{M_0}\varphi <\sup\{\sup_{M_t}\varphi: 0<t\leq T\}=\varphi
(x_0).
\end{equation*}
We choose a local normal coordinate system $(\xi ^i)$ at the point
$x_0\in M(t_0)$ such that at $x_0=x(t_0, \xi _0)$ we have
$g_{ij}=\delta _{ij}$, $ \varphi =h_n^n$ and $h_{ij}$
is diagonalized. Consider the contravariant
vector field  $\widetilde{\eta }=(\widetilde{\eta }^i)=(0, 0, \cdots, 0, 1)$,
and define on a neighbourhood of $(t_0, \xi_0)$ the function
\begin{equation*}
\widetilde{\varphi}=\frac {h_{ij}\widetilde{\eta
}^i\widetilde{\eta }^j}{g_{ij}\widetilde{\eta} ^i\widetilde{\eta}
^j}.
\end{equation*}
$\widetilde {\varphi }$ assumes its maximum at
$(t_0, \xi_0)$, and at this point,
$\dot {\widetilde{\varphi }}=\dot h _n^n$
and the spatial derivatives do coincide. So,  at $(t_0,
\xi_0)$, $\widetilde {\varphi }$ satisfies the same differential
equation (3.4) as $h_n^n$.
From the estimates in preceding sections, we deduce
from Lemma 3.3 that
\begin{equation}
\dot {\varphi }-pH^{p-1}\Delta \varphi \leq
pH^{p-1}\left(c(1+h_n^n)+(h_n^n)^2-||A||^2h_n^n\right),
\end{equation}
where  $0<p\leq 1$ and again $H$ is bounded in
finite time. From (6.1)
and the maximum principle we deduce
that $\varphi$ is uniformly bounded.\\
\end{proof}
Suppose that $T^*<\infty$, in the flow equation (4.3).
 Then from Lemma 4.1, we know that the
 flow stays in a compact region of $N$, i.e $u$ stays uniformly bounded.
Furthermore, in view of  Theorem 5.1 the first derivative  $Du$
stays also  uniformly bounded. Finally
Lemma 6.1 (with (2.8), (2.9) and (2.10)) we obtain uniform
 $C^2$-estimate for $u$.
Since in finite time, $H$ is bounded from above by (4.4)
 and (4.5), and we also know that $H$ has a positive lower
bound by Lemma 3.4, applying the standard
 regularity results \cite {k} (see also \cite {z}) we obtain
$C^{2, \alpha}$-estimates from
 the $C^2$-estimates of $u$, leading further to uniform
$C^{m,\alpha}$-estimates
 for any positive integer $m$.
But this contradicts to the maximality of $T^*$.
 Therefore, $T^*=\infty$, i.e., the flow exists for all time.

 Now we further assume that the
limit $\lim_{t\rightarrow \infty}\inf_{x\in S_0}u(t,x)$ is bounded
from below. Then from the proof of Lemma 4.1 we know that the flow
hypersurfaces $M(t)$ stay in a precompact region of $N$ for all
$t$. Integrating (4.3) with respect to $t$, and observing that the
right hand side is non positive, yields
$$u(0, x)-u(t, x)=\int _0^te^{-\psi}v(H^p-\tau)\geq c\int
_0^t(H^p-\tau),$$ i. e.
\begin{equation*}
\int_0^{\infty}(H^p-\tau)<\infty, \quad \forall x\in S_0.
\end{equation*}
Hence for any $x\in S_0$, there is a sequence $t_k\rightarrow
\infty $ such that, $H^p \rightarrow \tau$. On the other hand,
$u(\cdot, x )$ is monotone decreasing and therefore
\begin{equation*}
\lim _{t\rightarrow \infty}u(t,x)=\widetilde{u}(x)
\end{equation*}
exists and has desired regularity. We therefore have proved the
following theorem
\begin{theorem}
Let $N$ be a cosmological spacetime with compact Cauchy
hypersurface $S_0$ and $M_0$ a closed spacelike hypersurface in
$N$ such that $H\geq \sqrt[p]{\tau}$. If the curvature condition
(1.2) is satisfied, then the evolution equation (3.1) with initial
hypersurface $M_0$ exists for all time, when $p\in(0,1]$. Moreover
if we further assume that the limit $\lim_{t\rightarrow
\infty}\inf_{x\in S_0}u(t,x)$ is bounded from below,  then the
solutions to (3.1) converge to a stationary spacelike hypersurface
with mean curvature $\sqrt[p]{\tau}$, as $t\rightarrow \infty$.
\end{theorem}

The positive number $\tau$ is arbitrary, letting $\tau \rightarrow
0$ in the above theorem, we obtain that the power mean curvature
flow (1.1) exists for all time. If we further assume the limit
$\lim_{t\rightarrow \infty}\inf_{x\in S_0}u(t,x)$ is bounded from
below, the flow hypersurfces $M(t)$ of (1.1) converge to a
stationary maximum spacelike hypersurface. This completes the
proof of Theorem 1.1.

 \end{document}